\documentclass[12pt,twoside]{article}

\usepackage{headerfo}
\usepackage{amsfonts}

\newcounter{theorem}
\newcommand{\newsection}[1]{{\setcounter{theorem}{0} \section{#1}}}
\newtheorem{Theorem}{Theorem}[section]
\newtheorem{Definition}[Theorem]{Definition}
\newtheorem{Proposition}[Theorem]{Proposition}
\newtheorem{Lemma}[Theorem]{Lemma}
\newtheorem{Corollary}[Theorem]{Corollary}

\def\eop{{ \vrule height7pt width7pt depth0pt}\par\bigskip}

\newcommand{\R}{\mathbb R}
\newcommand{\C}{\mathbb C}

\newif\ifpdf
\ifx\pdfoutput\undefined
  \pdffalse
\else
  \pdfoutput=1
  \pdftrue
\fi

\usepackage{latexsym}

\chardef\aa=64

\begin{document}

\renewcommand{\author} {L.~Bos \\  
Department of Computer Science\\
University of Verona \\
Verona, Italy\\
\medskip
N.~Levenberg\\
Department of Mathematics \\
Indiana University \\
Bloomington, Indiana, \\
USA\\
\medskip
and \\
J.~Ortega-Cerd\`a \\
Department of Mathematics \\
University of Barcelona \\
Barcelona, Spain\\}

\renewcommand{\title}{Optimal Polynomial Prediction Measures and Extremal Polynomial Growth}
\newcommand{\stitle}{Optimal Prediction}

\renewcommand{\date}{\today}
\flushbottom
\setcounter{page}{1}
\pageheaderlinetrue
\oddpageheader{}{\stitle}{\thepage}
\evenpageheader{\thepage}{\stitle}{}
\thispagestyle{empty}
\vskip1cm
\begin{center}
\LARGE{\bf \title}
\\[0.7cm]
\large{\author}
\end{center}
\vspace{0.3cm}
\begin{center}
\date
\end{center}
\begin{abstract}
We show that the problem of finding the measure supported on a compact set $K\subset \C$ such that the 
variance of the least squares predictor by polynomials of degree at most $n$ at a point $z_0\in\C^d\backslash K$ is a minimum, is equivalent to the problem of finding the polynomial of degree at most $n,$ bounded by 1 on $K,$ with extremal growth at $z_0.$ We use this to find the polynomials of extremal growth for $[-1,1]\subset \C$ at a purely imaginary point. The related  problem on the extremal growth of real polynomials was studied by Erd\H{o}s in 1947, \cite{E}.
\end{abstract}

\setcounter{section}{0}

\newsection{Introduction} In this work we consider two classical extremum problems for polynomials. The first is very easy to state. Indeed, let us denote the {\it complex} polynomials of degree at most $n$ in $d$ complex variables 
 by $\C_n[z],$ $z\in \C^d.$ Then for $K\subset \C^d$ compact and $z_0\in\C^d\backslash K$ an {\it external}  point, we say that $P_n(z)\in \C_n[z]$ has {\it extremal growth relative to $K$} at $z_0$ if the sup-norm $\|P_n\|_K\le1$ and
\begin{equation}\label{ExtGrowthDef}
|P_n(z_0)|=\max_{p\in \C_n[z],\,\,\|p\|_K\le1}|p(z_0)|.
\end{equation}
We note that for this to be well-defined we require that $K$ be polynomial determining, i.e., if $p\in\C[z]$ is such that $p(x)=0$ for all $x\in K,$ then $p=0.$ We refer the interested reader to the survey \cite{BMW} for more about what is known about this problem.

\medskip

The second problem is from the field of Optimal Design for Polynomial Regression. To describe it we reduce to the real case $K\subset \R^d,$ and note that we may write any $p\in \R_n[z]$ in the form
\[p=\sum_{k=1}^N \theta_k p_k\]
where ${\cal B}_s:=\{p_1,p_2,\ldots,p_N\}$ is a basis for $\R_n[z]$ and $N:={n+d\choose d}$ its dimension.

Suppose now that we observe the values of a particular
$p\in \R_n[z]$ at a set of  $m\ge N$ points $X:=\{x_j\,:\, 1\le j\le m\}\subset K$ with some random errors, i.e., we observe
\[y_j =p(x_j)+\epsilon_j,\quad 1\le j\le m\]
where we assume that the errors $\epsilon_j\sim N(0,\sigma)$ are independent. In matrix form this becomes
\[ y= V_n\theta+\epsilon\]
where $\theta \in\R^N,$ $y,\epsilon \in \R^m$ and
\[V_n:=\left[\begin{array}{cccccc}
p_1(x_1)&p_2(x_1)&\cdot&\cdot&\cdot&p_N(x_1) \cr
p_1(x_2)&p_2(x_2)&\cdot&\cdot&\cdot&p_N(x_2) \cr
\cdot&&&&&\cdot\cr
\cdot&&&&&\cdot\cr
\cdot&&&&&\cdot\cr
\cdot&&&&&\cdot\cr
\cdot&&&&&\cdot\cr
p_1(x_m)&p_2(x_m)&\cdot&\cdot&\cdot&p_N(x_m) \end{array}\right]\in \R^{m\times N}\]
is the associated Vandermonde matrix.

Our assumption on the error vector $\epsilon$ means that
\[{\rm cov}(\epsilon)=\sigma^2I_m\in\R^{m\times m}.\]
Now, the least squares estimate of $\theta$ is
\[\widehat{\theta}:=(V_n^tV_n)^{-1}V_n^ty.\]

Note  that the entries of $\displaystyle{ {1\over m}V_n^tV_n}$ are the discrete inner products of the $p_i$
with respect to the measure
\begin{equation}\label{StatsMeas}
\mu_X ={1\over m}\sum_{k=1}^m \delta_{x_k}.
\end{equation}
More specifically,
\[ {1\over m}V_n^tV_n=G_n(\mu_X)\]
where
\begin{equation}\label{GReal}
G_n(\mu):=\left[\int_K p_i(x)p_j(x)d\mu\right]_{1\le i,j\le N}\in\R^{N\times N}
\end{equation}
is the  Moment, or Gram, matrix of the polynomials $p_i$ with respect to the measure $\mu.$

In general we may consider arbitrary probability measures on $K,$ setting
\[{\cal M}(K):=\{\mu\,:\, \mu \,\,\hbox{is a probability measure on}\,\, K\}.\]

Now set
\begin{equation}\label{P}
{\mathbf p}(z)=\left[\begin{array}{c}p_1(z)\cr p_2(z)\cr\cdot\cr\cdot\cr p_N(z)\end{array}\right]\in \R^N
\end{equation}
then the least squares estimate of the {\it observed polynomial} is
\[{\mathbf p}^t(z)\widehat{\theta}.\]
We may compute its variance at any point $z\in \R^d$ to be
\begin{eqnarray}
{\rm var}({\mathbf p}^t(z)\widehat{\theta})&=&\sigma^2{\mathbf p}^t(z)(V_n^tV_n)^{-1}{\mathbf p}(z) \nonumber\\
&=&{1\over m}\sigma^2 {\mathbf p}^t(z)(G_n(\mu_X))^{-1}{\mathbf p}(z)\label{variance}
\end{eqnarray}
where $\mu_X$ is again given by (\ref{StatsMeas}). Now, it is easy to verify that for any $\mu\in{\cal M}(K),$
\[ {\mathbf p}^t(z)(G_n(\mu))^{-1}{\mathbf p}(z)=K_n^\mu(z,z) \]
where, for $\{q_1,\cdots,q_N\}\subset \R_n[z],$ a $\mu$-orthonormal basis for $\R_n[z],$
\[ K_n^\mu(w,z):=\sum_{k=1}^N \overline{q_k(w)}q_k(z) \]
is the Bergman kernel for $\R_n[z]$. The function $K_n^\mu(z,z)$ is also known as the 
(reciprocal of) the Christoffel function for $\R_n[z].$

\medskip
We may generalize easily to the complex case, $K\subset \C^d,$ where now the $p_j$ form a basis for
$\C_n[z]$ and
\begin{equation}\label{GComplex}
G_n(\mu):=\left[\int_K p_i(z)\overline{p_j(z)}d\mu\right]_{1\le i,j\le N}\in\C^{N\times N}.
\end{equation}
For an external point $z_0\in \C^d\backslash K,$ a measure $\mu_{0}\in {\cal M}(K)$ is said to be an {\it optimal prediction (or extrapolation) measure for $z_0$ relative to $K$} if it minimizes the complex analogue of the variance (\ref{variance}) of the polynomial predictor at $z_0,$ i.e., if
\begin{equation}\label{OptMeasDef}
K_n^{\mu_0}(z_0,z_0)= \min_{\mu \in {\cal M}(K)}
K_n^{\mu}(z_0,z_0).
\end{equation}

\bigskip
In \cite{HL} Hoel and Levine show that in the univariate case, for $K=[-1,1],$ and {\it any} $z_0\in \R\backslash K,$  a {\it real} external point, the optimal prediction measure is a discrete measure supported at the $n+1$ extremal points $x_k=\cos(k\pi/n),$ $0\le k\le n,$ of $T_n(x)$ the classical Chebyshev polynomial of the first kind (cf. Lemma \ref{HLWeights} below). In this case it turns out that
\begin{equation}\label{HLthm}
K_n^{\mu_0}(z_0,z_0)=T_n^2(z_0).
\end{equation}
Notably, as is well known, $T_n(x)$ is the polynomial of extremal growth for any point $z_0\in\R\backslash [-1,1]$ relative to $K=[-1,1].$ Also, Erd\H{o}s (1947) \cite{E} has shown that the Chebyshev polynomial is also extreme relative to $[-1,1]$ for {\it real} polynomials at points $z_0\in\C$ with $|z_0|\ge1,$ i.e.,
\[\max_{p\in\R_n[x],\,\,\|p\|_{[-1,1]}\le1}|p(z_0)|=|T_n(z_0)|.\] The problem for real polynomials and $|z_0|\le1$ or for complex polynomials $p\in\C[z]$ has remained unsolved up to now.

We show in Section \ref{Equivalence} that (\ref{HLthm}) is not an accident, and that there is a general equivalence of our two extremum problems. In Section \ref{ComplexPoint} we will use this to compute the polynomials of extremal growth and the optimal prediction measures for a purely imaginary complex point $z_0\in\C\backslash [-1,1].$

\section{A Kiefer-Wolfowitz Type Equivalence Theorem}\label{Equivalence}

Kiefer and Wolfowitz \cite{KW} have given a remarkable equivalence between what are called D-optimal and G-optimal designs, i.e., probability measures that maximize the determinant of the {\it design} matrix $G_n(\mu)$ and those which minimize the maximum over $x$ {\it interior} to $K,$ of the prediction variance
  i.e., minimize $\max_{x\in K} K_n^\mu(x,x).$ Here we give an analogous equivalence, for a single 
{\it exterior} point $z_0\in \C^d\backslash K,$ with the problem of extremal polynomial growth.

Of importance will be the well-known variational form of the Christoffel function:
\begin{equation}\label{VarForm}
K_n^\mu(z_0,z_0)=\sup_{p\in \C_n[z]}\frac{|p(z_0)|^2}{\int_K |p(z)|^2d\mu}.
\end{equation}

Indeed, from this variational form, the problem of minimal variance
(\ref{OptMeasDef})
may be expressed as
\[ \min_{\mu \in {\cal M}(K)} \max_{p\in \C_n[z]}\frac{|p(z_0)|^2}{\int_K |p(z)|^2d\mu}\]
which, as it turns out, we will be able to analyze using the classical Minimax Theorem (see e.g. Gamelin \cite[Thm. 7.1, Ch. II] {G}). To see this,
note that we may first of all simplify to
\[
\min_{\mu \in {\cal M}(K)} \max_{p\in \C_n[z],\,p(z_0)=1}\frac{1}{\int_K |p(z)|^2d\mu}
=1/\left\{\max_{\mu \in {\cal M}(K)} \min_{p\in \C_n[z],\,p(z_0)=1}\int_K |p(z)|^2d\mu\right\}.
\]
Now, for $\mu\in {\cal M}(K)$ and $p\in \C_n[z]$ such that $p(z_0)=1,$ let
\[f(\mu,p):=\int_K |p(z)|^2d\mu.\]
It is easy to confirm that $f$ is quasiconcave in $\mu$ and quasiconvex in $p$ and hence by the Minimax Theorem
\[\max_{\mu \in {\cal M}(K)} \min_{p\in \C_n[z],\,p(z_0)=1}\int_K |p(z)|^2d\mu= \min_{p\in \C_n[z],\,p(z_0)=1} \max_{\mu \in {\cal M}(K)} \int_K |p(z)|^2d\mu.\]
However, as $\mu=\delta_x\in {\cal M}(K)$ for every $x\in K,$ it
follows that 
\[\max_{\mu \in {\cal M}(K)} \int_K |p(z)|^2d\mu=\|p\|_K^2.\]
Consequently, the minimum variance is given by
\[ \min_{\mu \in {\cal M}(K)} \max_{p\in \C_n[z]}\frac{|p(z_0)|^2}{\int_K |p(z)|^2d\mu}=\max_{p\in \C_n[z]} \frac{|p(z_0)|^2}
{\|p\|_K^2},\]
i.e., the value squared of the polynomial of maximal growth at $z_0$. The Minimax theorem in a similar context has been used before to get pointwise estimates of solutions to the $\bar\partial$-equation by Berndtsson in \cite[p. 206]{Berndtsson}. 
\medskip

It is also possible to give a more precise relation between the
extremal polynomials for the two problems (of minimum variance and extremal growth), using completely elementary means, along the lines of the proof of the original Kiefer-Wolfowitz theorem.

We begin with a simple technical lemma.

\begin{Lemma}\label{NormCond}
Suppose that $\mu\in{\cal M}(K)$ is such that $G_n(\mu)$ is non-singular and $z_0\in\C^d\backslash K.$ Define the polynomial
\begin{equation}\label{DefP}
P_n^{\mu,z_0}(z):=\frac{{\mathbf p}^*(z_0)G_n^{-1}(\mu){\mathbf p}(z)}
{\sqrt{{\mathbf p}^*(z_0)G_n^{-1}(\mu){\mathbf p}(z_0)}}
=\frac{K_n^\mu(z_0,z)}{\sqrt{K_n^\mu(z_0,z_0)}}.
\end{equation}
Then if $\|P_n^{\mu,z_0}\|_K\le1$ it is a polynomial of degree at most $n$ of extremal growth at $z_0$ relative to $K.$ Here ${\mathbf p}^*$ denotes the conjugate transpose of the vector ${\mathbf p}.$
\end{Lemma}
\noindent {\bf Proof.} Suppose that $p\in\C_n[z]$ is such that $\|p\|_K\le1.$ Then, from the variational form
(\ref{VarForm}), 
\begin{eqnarray*}
|p(z_0)|^2&\le &K_n^\mu(z_0,z_0)\times \int_K |p(z)|^2d\mu\cr
&\le&K_n^\mu(z_0,z_0)\times \int_K \|p\|_K^2 d\mu\cr
&\le&K_n^\mu(z_0,z_0)\cr
&=&|P_n^{\mu,z_0}(z_0)|^2. 
\end{eqnarray*}
\noindent \eop
\begin{Theorem} \label{EquivThm}
Suppose that $\mu_0\in{\cal M}(K)$ and define $P_n^{\mu_0,z_0}$ as in (\ref{DefP}). Then
$\mu_0$ is an extremal prediction measure for $z_0$ relative to $K$ if and only if $P_n^{\mu_0,z_0}(z)\in \C_n[z]$
is a polynomial of extremal growth at $z_0$ relative to $K$ (i.e., $\|P_n^{\mu_0,z_0}\|_K\le1$, by Lemma \ref{NormCond}).
\end{Theorem}
\noindent {\bf Proof.} Suppose first that $\|P_n^{\mu_0,z_0}\|_K\le1.$ We will show that then $\mu_0$ is an optimal prediction measure. Indeed, consider any other measure $\mu\in{\cal M}(K).$ Then
\begin{eqnarray*}
K_n^\mu(z_0,z_0)&=& \sup_{p\in \C_n[z]}\frac{|p(z_0)|^2}{\int_K |p(z)|^2d\mu} \cr
&\ge& \frac{|P_n^{\mu_0,z_0}(z_0)|^2}{\int_K |P_n^{\mu_0,z_0}(z)|^2d\mu}\cr
&\ge& \frac{|P_n^{\mu_0,z_0}(z_0)|^2}{\int_K \|P_n^{\mu_0,z_0}\|_K^2d\mu}\cr
&\ge& \frac{|P_n^{\mu_0,z_0}(z_0)|^2}{\int_K 1\,d\mu}\cr
&=&|P_n^{\mu_0,z_0}(z_0)|^2\cr
&=&K_n^{\mu_0}(z_0,z_0).
\end{eqnarray*}
Conversely, suppose that $\mu_0$ is an extremal prediction measure. We must show that then
$\|P_n^{\mu_0,z_0}\|_K\le 1.$ To see this, fix $\mu_1\in {\cal M}(K)$ and consider the family of
measures $d\mu_t:=t\,d\mu_1+(1-t)\,d\mu_0,$ $0\le t\le 1.$ We calculate
\begin{eqnarray*}
\frac{d}{dt} K_n^{\mu_t}(z_0,z_0)&=&\frac{d}{dt}\left\{{\mathbf p}^*(z_0)G_n^{-1}(\mu_t){\mathbf p}(z_0)\right\}\cr
&=&{\mathbf p}^*(z_0)\left\{\frac{d}{dt}G_n^{-1}(\mu_t)\right\}{\mathbf p}(z_0)\cr
&=&{\mathbf p}^*(z_0)\left\{-G_n^{-1}(\mu_t)\left(\frac{d}{dt}G_n(\mu_t)\right)
G_n^{-1}(\mu_t)\right\}{\mathbf p}(z_0).
\end{eqnarray*}
But,
\begin{eqnarray*}
\frac{d}{dt}G_n(\mu_t)&=&\frac{d}{dt} [\int_K \overline{p_i(z)}p_j(z) (td\mu_1+(1-t)d\mu_0)]\cr
&=& [\int_K \overline{p_i(z)}p_j(z) (d\mu_1-d\mu_0)]\cr
&=&G_n(\mu_1)-G_n(\mu_0).
\end{eqnarray*}
Hence,
\begin{eqnarray*}
\left.\frac{d}{dt} K_n^{\mu_t}(z_0,z_0)\right|_{t=0}
&=&{\mathbf p}^*(z_0)\left\{-G_n^{-1}(\mu_0)\left(G_n(\mu_1)-G_n(\mu_0)\right) \times \right.\cr
&& \qquad \left.G_n^{-1}(\mu_0)\right\}{\mathbf p}(z_0)\cr
&=&{\mathbf p}^*(z_0)G_n^{-1}(\mu_0){\mathbf p}(z_0) \cr 
&&- {\mathbf p}^*(z_0)G_n^{-1}(\mu_0) G_n(\mu_1)G_n^{-1}(\mu_0){\mathbf p}(z_0)\cr
&=&K_n^{\mu_0}(z_0,z_0) \cr
&&- {\mathbf p}^*(z_0)G_n^{-1}(\mu_0) G_n(\mu_1)G_n^{-1}(\mu_0){\mathbf p}(z_0).
\end{eqnarray*}
Now, for a fixed $a\in K,$ take $\mu_1=\delta_a,$ the Dirac delta measure supported at $a.$ In this case
\[G_n(\mu_1)=[p_i(a)\overline{p_j(a)}]\in \C^{N\times N}={\mathbf p}(a){\mathbf p}^*(a).\]
Hence, using (\ref{DefP}),
\begin{eqnarray*}
\left.\frac{d}{dt} K_n^{\mu_t}(z_0,z_0)\right|_{t=0}&=&K_n^{\mu_0}(z_0,z_0)
-{\mathbf p}^*(z_0)G_n^{-1}(\mu_0) ({\mathbf p}(a){\mathbf p}^*(a))
G_n^{-1}(\mu_0){\mathbf p}(z_0)\cr
&=&K_n^{\mu_0}(z_0,z_0)-\left({\mathbf p}^*(z_0)G_n^{-1}(\mu_0) {\mathbf p}(a)\right)\left(
{\mathbf p}^*(a)G_n^{-1}(\mu_0){\mathbf p}(z_0)\right)\cr
&=&K_n^{\mu_0}(z_0,z_0)-P_n^{\mu_0,z_0}(a)\overline{P_n^{\mu_0,z_0}(a)}K_n^{\mu_0}(z_0,z_0)\cr
&=&K_n^{\mu_0}(z_0,z_0)\left(1-|P_n^{\mu_0,z_0}(a)|^2\right).
\end{eqnarray*}
If $\mu_0$ is to minimize $K_n^{\mu_t}(z_0,z_0)$ then each of the above derivatives must be greater than or equal to zero, i.e.,
\[|P_n^{\mu_0,z_0}(a)|\le1,\quad \forall a\in K.\]
\noindent \eop

\bigskip
\noindent {\bf Remark.} It is easily confirmed that $\int_K |P_n^{\mu,z_0}(z)|^2d\mu(z)=1.$ Hence, for an optimal prediction measure $\mu_0,$ it must be the case that 
$|P_n^{\mu_0,z_0}(z)|\equiv 1$ on the support of $\mu_0.$ Consequently optimal prediction measures are always supported on a real algebraic subset of $K$ of degree $2n.$ \eop

\section{A Complex Point External to $[-1,1]$}\label{ComplexPoint}

We now consider $K=[-1,1]\subset\C$ and $z_0\in\C\backslash K.$ First note that by the above remark any optimal prediction measure must be supported on discrete points, $x_0:=-1,$ $x_n:=+1$ and $n-1$ internal points
$-1<x_1<\cdots<x_{n-1}<1,$ i.e., is of the form 
\[\mu_0=\sum_{i=0}^n w_i \delta_{x_i}\]
with weights $w_i>0,$ $\sum_{i=0}^n w_i=1.$ Given the support points $x_i$ there is a simple recipe for the optimal weights, given already in \cite{HL}.

\begin{Lemma} \label{HLWeights} (Hoel-Levine) Suppose that $-1=x_0<x_1<\cdots<x_n=+1$ are given. Then among all discrete probability measures supported at these points, the measure with
\begin{equation}\label{OptWeights}
w_i:=\frac{|\ell_i(z_0)|}{\sum_{i=0}^n |\ell_i(z_0)|},\,\,0\le i\le n
\end{equation}
with $\ell_i(z)$ the $i$th fundamental Lagrange interpolating polynomial for these points, minimizes
$K_n^\mu(z_0,z_0).$
\end{Lemma}
\noindent {\bf Proof.} 
We first note that for such a discrete measure, $\{\ell_i(z)/\sqrt{w_i}\}_{0\le i \le n}$ form an orthonormal basis. Hence
\begin{equation}\label{GenWeights}
K_n^\mu(z_0,z_0)=\sum_{i=0}^n \frac{|\ell_i(z_0)|^2}{w_i}.
\end{equation}
In the case of the weights chosen according to (\ref{OptWeights}) we obtain
\begin{equation}\label{OptKn}
K_n^{\mu_0}(z_0,z_0)=\left(\sum_{i=0}^n |\ell_i(z_0)|\right)^2.
\end{equation}
We claim that for any choice of weights $K_n$ given by (\ref{GenWeights}) is at least as large as that given by (\ref{OptKn}). To see this, just note that by the Cauchy-Schwartz inequality,
\begin{eqnarray*}
\left(\sum_{i=0}^n |\ell_i(z_0)|\right)^2&=&\left(\sum_{i=0}^n \frac{|\ell_i(z_0)|}{\sqrt{w_i}}\cdot
\sqrt{w_i}\right)^2\cr
&\le& \left(\sum_{i=0}^n \frac{|\ell_i(z_0)|^2}{w_i}\right)\cdot\left(\sum_{i=0}^n w_i\right)\cr
&=&\sum_{i=0}^n \frac{|\ell_i(z_0)|^2}{w_i}.
\end{eqnarray*}
\noindent \eop

\bigskip
\noindent {\bf Remark.} We note that the optimal $K_n(z_0,z_0)$ given by (\ref{OptKn}) is the Lebesgue function squared. Hence the problem of finding the support of the optimal prediction measure amounts to finding the $n+1$ interpolation points $-1=x_0<x_1<\cdots<x_n=+1$ for which the Lebesgue function evaluated at the external point $z_0,$
\[\Lambda(z_0):=\sum_{i=0}^n |\ell_i(z_0)|,\]
is as small as possible. \eop

\medskip

In this case the extremal polynomial $P_n^{\mu,z_0}(z)$ also simplifies.

\begin{Lemma}\label{SimpleP} Suppose that the measure $\mu_0$ is supported at the points
$-1=x_0<x_1<\cdots<x_n=+1$ with optimal weights given by (\ref{OptWeights}). Then
\[P_n^{\mu_0,z_0}(z)=\sum_{i=0}^n {\rm sgn}(\ell_i(z_0))\ell_i(z)\]
where ${\rm sgn}(z):=\overline{z}/|z|$ is the complex sign of $z\in\C.$
\end{Lemma}
\noindent {\bf Proof.} Using again the fact that $\{\ell_i(z)/\sqrt{w_i}\}_{0\le i\le n}$ form a set of
orthonormal polynomials, we have
\begin{eqnarray*}
P_n^{\mu_0,z_0}(z)&=&\frac{1}{\Lambda(z_0)}\sum_{i=0}^n \frac{\overline{\ell_i(z_0)}}{\sqrt{w_i}}\frac{\ell_i(z)}{\sqrt{w_i}}\cr
&=&\frac{1}{\Lambda(z_0)}\sum_{i=0}^n \left(\Lambda(z_0)  \frac{\overline{\ell_i(z_0)}}{|\ell_i(z_0)|}
\right)\ell_i(z)\cr
&=& \sum_{i=0}^n \frac{\overline{\ell_i(z_0)}}{|\ell_i(z_0)|} \cdot \ell_i(z).
\end{eqnarray*}
\noindent \eop

\medskip
\noindent {\bf Remark.} By the equivalence Theorem \ref{EquivThm} the support of the optimal prediction measure and the polynomial of extremal growth will be given by those points $-1=x_0<x_1<\cdots<x_n=+1$ for which 
\[\max_{-1\le x\le 1} \left|\sum_{i=0}^n \frac{\overline{\ell_i(z_0)}}{|\ell_i(z_0)|} \cdot \ell_i(x)\right|=1.\]
\noindent \eop

\section{A Purely Imaginary Point External to $[-1,1]$}
In the case of $z_0=ai,$ $0\neq a\in\R,$ a purely imaginary point, it turns out that there are remarkable formulas for the polynomial of extremal growth as well as for the support of the optimal prediction measure. Both of these will depend on the point $z_0$ (as opposed to the real case $z_0\in \R\backslash[-1,1]$ where Hoel and Levine \cite{HL} showed that the support is always the set of extreme points of the Chebyshev polynomial $T_n(x)$).

To begin we will first analyze the degrees $n=1$ and $n=2$ cases.

\subsection{Degree $n=1$}
Here the support of the extremal measure is necessarily $x=-1$ and $x_1=+1.$ We will compute 
$P_1^{\mu_0,z_0}(z)$ using the formula given in Lemma \ref{SimpleP}. Indeed in this case,
$\ell_0(z)=(1-z)/2$ and $\ell_1(z)=(1+z)/2$ so that
\[{\rm sgn}(\ell_0(ia))={\rm sgn}\left(\frac{1-ia}{2}\right)=\frac{1+ia}{\sqrt{a^2+1}}\]
and
\[{\rm sgn}(\ell_1(ia))={\rm sgn}\left(\frac{1+ia}{2}\right)=\frac{1-ia}{\sqrt{a^2+1}}.\]
Hence,
\begin{eqnarray*}
P_1^{\mu_0,z_0}(z)&=& \frac{1+ia}{\sqrt{a^2+1}}\frac{1-z}{2}
+\frac{1-ia}{\sqrt{a^2+1}}\frac{1+z}{2}\cr
&=&\frac{1}{\sqrt{a^2+1}}\{1-iaz\}.
\end{eqnarray*}
Since $\pm1$ is necessarily the support of the optimal prediction measure it is immediate that 
$\|P_1^{\mu_0,z_0}\|_{[-1,1]}=1,$ as is also easily verified by a simple direct calculation.

\subsection{Degree $n=2$}
We claim that the support of the optimal prediction measure is $x_0=-1,$ $x_1=0$ and $x_2=+1.$ However, this is not automatic and we will have to verify that the norm of $P_2^{\mu_0,z_0}$ is indeed 1.
Now, it is easy to see, for this support,  that
\[\ell_0(z)=\frac{z(z-1)}{2},\quad \ell_1(z)=1-z^2,\quad \ell_2(z)=\frac{z(z+1)}{2}\]
for which
\begin{eqnarray*}
{\rm sgn}(\ell_0(ia))&=&{\rm sgn}\left(\frac{ia(ia-1)}{2}\right)\cr
&=&\frac{-ia}{|a|}\cdot \frac{-ia-1}{\sqrt{a^2+1}}\cr
&=& i\,{\rm sgn}(a)\frac{1+ia}{\sqrt{a^2+1}},
\end{eqnarray*}
\[{\rm sgn}(\ell_1(ia))={\rm sgn}(1+a^2)=+1,\]
and, after a simple calculation,
\[{\rm sgn}(\ell_2(ia))= i\,{\rm sgn}(a)\frac{ia-1}{\sqrt{a^2+1}}.\]
From this we may easily conclude that
\begin{eqnarray*}
P_2^{\mu_0,z_0}(z)&=&\sum_{i=0}^2 {\rm sgn}(\ell_i(z_0))\ell_i(z)\cr
&=& \frac{{\rm sgn}(a)}{\sqrt{a^2+1}}\left(
-(a+{\rm sgn}(a)\sqrt{a^2+1})z^2-iz+{\rm sgn}(a)\sqrt{a^2+1}\right).
\end{eqnarray*}
The fact that $\|P_2^{\mu_0,z_0}\|_{[-1,1]}=1$  is an immediate consequence of the following lemma.

\begin{Lemma} For $x\in\R$ we have
\begin{eqnarray*}
|P_2^{\mu_0,z_0}(x)|^2&=&1+\frac{(|a|+\sqrt{a^2+1})^2}{a^2+1}x^2(x^2-1)\cr
&=& 1+(x^2-1)R_1^2(x),\quad R_1(x):=\frac{|a|+\sqrt{a^2+1}}{\sqrt{a^2+1}}\,x.
\end{eqnarray*}
\end{Lemma}
\noindent {\bf Proof.} This follows from elementary calculations starting with the formula for $P_2^{\mu_0,z_0}(x)$ given above. \eop

\medskip

We now define a sequence of polynomials, $Q_n(z),$ based on the above degrees $n=1$ and $n=2$ cases,
for which we will show that $Q_n(z)=c_nP_n^{\mu_0,z_0}(z)$ for certain $c_n\in\C$ with modulus $|c_n|=1.$ We will also define a sequence of polynomials $R_n(x)$ which will play the role of $R_1(x)$ in the Lemma for general degree $n.$

Now, as the formula for $P_2^{\mu_0,z_0}$ depends on the sign of $a,$ in order to simplify the formulas we will assume that $a>0.$ For $a<0,$ one may use the relation $P_2^{\mu_0,ia}(z)=P_2^{\mu_0,-ia}(-z).$

\begin{Definition} For $a>0$ we define the sequences of polynomials $Q_n(z)$ and $R_n(z)$ by
\begin{eqnarray*}
Q_1(z)&=& -\frac{az+i}{\sqrt{a^2+1}}, \qquad (=(-i)P_1^{\mu_0,z_0}(z))\cr
Q_2(z)&=&\frac{1}{\sqrt{a^2+1}}\left(-(a+\sqrt{a^2+1})z^2-iz+\sqrt{a^2+1}\right),
\quad (=P_2^{\mu_0,z_0}(z))\cr
Q_{n+1}(z)&=& 2zQ_n(z)-Q_{n-1}(z),\quad n=2,3,\cdots.
\end{eqnarray*}
and
\begin{eqnarray*}
R_0(z)&=&\frac{a}{\sqrt{a^2+1}},\cr
R_1(z)&=&\frac{a+\sqrt{a^2+1}}{\sqrt{a^2+1}} \,z,\cr
R_{n+1}(z)&=&2zR_n(z)-R_{n-1}(z),\quad n=1,2,\cdots.
\end{eqnarray*}
\end{Definition}
Since the recursions are both those of the classical Chebyshev polynomials it is not surprising that there
are formulas for $Q_n(z)$ and $R_n(z)$ in terms of these.

\begin{Lemma}\label{Qn} We have
\[
Q_n(z)=\frac{1}{\sqrt{a^2+1}}\left( -(az+i)T_{n-1}(z)+\sqrt{a^2+1}(1-z^2)U_{n-2}(z)\right)
\]
where $T_n(z)$ is Chebyshev polynomial of the first kind and $U_n(z):=\frac{1}{n+1}T_{n+1}'(z)$ that of the second kind.
\end{Lemma}
\noindent {\bf Proof.} Let $q_n(z)$ denote the right side of the proposed identity. We proceed by induction.
For $n=1$ we have 
\begin{eqnarray*}
q_1(z)&=&\frac{1}{\sqrt{a^2+1}}\left( -(az+i)T_{1-1}(z)+\sqrt{a^2+1}(1-z^2)U_{1-2}(z)\right)\cr
&=&\frac{1}{\sqrt{a^2+1}}\left(-(az+i)\times1 + 0\right)\cr
&=&Q_1(z).
\end{eqnarray*}
Similarly, for $n=2$ we have
\begin{eqnarray*}
q_2(z)&=&\frac{1}{\sqrt{a^2+1}}\left( -(az+i)T_{2-1}(z)+\sqrt{a^2+1}(1-z^2)U_{2-2}(z)\right)\cr
&=&\frac{1}{\sqrt{a^2+1}}\left(-(az+i)z+\sqrt{a^2+1}(1-z^2)\right)\cr
&=&\frac{1}{\sqrt{a^2+1}}\left( -(a+\sqrt{a^2+1})z^2-iz+\sqrt{a^2+1} \right)\cr
&=&Q_2(z).
\end{eqnarray*}
The result now follows easily from the fact that both kinds of Chebyshev polynomials satisfy the same recursion as used in the definition of $Q_n(z).$ \eop
\medskip
\begin{Lemma}\label{Rn}
We have
\[R_n(z)=\frac{1}{\sqrt{a^2+1}} \left(\sqrt{a^2+1}zU_{n-1}(z)+aT_{n}(z)\right).\]
\end{Lemma}
\noindent {\bf Proof.} let $r_n(z)$ denote the right side of the proposed identity. We again proceed by induction. For $n=0$ we have
\begin{eqnarray*}
r_0(z)&=&\frac{1}{\sqrt{a^2+1}} \left(\sqrt{a^2+1}zU_{-1}(z)+aT_{0}(z)\right)\cr
&=&\frac{a}{\sqrt{a^2+1}}\cr
&=&R_0(z).
\end{eqnarray*}
Similarly, for $n=1$ we have
\begin{eqnarray*}
r_1(z)&=&\frac{1}{\sqrt{a^2+1}} \left(\sqrt{a^2+1}zU_{0}(z)+aT_{1}(z)\right)\cr
&=&\frac{1}{\sqrt{a^2+1}} \left(\sqrt{a^2+1}z\times1+a\times z\right)\cr
&=&\frac{a+\sqrt{a^2+1}}{\sqrt{a^2+1}}z\cr
&=&R_1(z).
\end{eqnarray*}
The result now follows easily from the fact that both kinds of Chebyshev polynomials satisfy the same recursion as used in the definition of $R_n(z).$ \eop
\medskip
Now, just for the Chebyshev polynomials $T_n(z)$ and $U_{n-1}(z)$ 
there is the Pell identity
\begin{equation}\label{ChebPell}
T_n^2(z)-(z^2-1)U_{n-1}^2(z)\equiv 1.
\end{equation}

We will show that for real $z\in\R,$ the polynomials
$Q_n(z)$ and $R_{n-1}(z)$ satisfy a similar Pell identity.

\begin{Proposition}\label{Pell} For $z=x\in\R,$ we have
\[|Q_n(x)|^2-(x^2-1)R_{n-1}^2(x)\equiv 1.\]
\end{Proposition}
\noindent {\bf Proof.} By Lemma \ref{Qn}, $z=x\in\R,$ we may write
\begin{eqnarray*}
Q_n(x)&=&\frac{1}{\sqrt{a^2+1}}\left( -(ax+i)T_{n-1}(x)+\sqrt{a^2+1}(1-x^2)U_{n-2}(x)\right)\cr
&=&\frac{1}{\sqrt{a^2+1}}\left(-iT_{n-1}(x)+\left\{-axT_{n-1}(x)+\sqrt{a^2+1}(1-x^2)U_{n-2}(x)
\right\}\right)\cr
\end{eqnarray*}
so that
\[
|Q_n(x)|^2=\frac{1}{a^2+1}\left(T_{n-1}^2(x)+\left(-axT_{n-1}(x)+\sqrt{a^2+1}(1-x^2)U_{n-2}(x)
\right)^2\right). \]
Hence, using the Chebyshev Pell identity (\ref{ChebPell}),
\begin{eqnarray*}
&(a^2+1)(1-|Q_n(x)|^2)&\cr
&=(a^2+1)-T_{n-1}^2(x)-a^2x^2T_{n-1}^2(x)&\cr
&-(a^2+1)(1-x^2)^2U_{n-2}^2(x)+2a\sqrt{a^2+1}x(1-x^2)U_{n-2}(x)T_{n-1}(x)&\cr
&=(a^2+1)(1-T_{n-1}^2(x))+a^2(1-x^2)T_{n-1}^2(x)-(a^2+1)(1-x^2)^2U_{n-2}^2(x)&\cr
&+2a\sqrt{a^2+1}x(1-x^2)U_{n-2}(x)T_{n-1}(x)&\cr
&=(a^2+1)(1-x^2)U_{n-2}^2(x)+a^2(1-x^2)T_{n-1}^2(x)-(a^2+1)(1-x^2)^2U_{n-2}^2(x)&\cr
&+2a\sqrt{a^2+1}x(1-x^2)U_{n-2}(x)T_{n-1}(x)&\cr
&=(1-x^2)\Big\{(a^2+1)U_{n-2}^2(x) +a^2T_{n-1}^2(x)-(a^2+1)(1-x^2)U_{n-2}^2(x)&\cr
&+2a\sqrt{a^2+1}xU_{n-2}(x)T_{n-1}(x)\Big\}&\cr
&=(1-x^2)\Big\{(a^2+1)[1-(1-x^2)]U_{n-2}^2(x) +a^2T_{n-1}^2(x)&\cr
&+2a\sqrt{a^2+1}xU_{n-2}(x)T_{n-1}(x)\Big\}&\cr
&=(1-x^2)\Big\{(a^2+1)x^2U_{n-2}^2(x) +a^2T_{n-1}^2(x)
+2a\sqrt{a^2+1}xU_{n-2}(x)T_{n-1}(x)\Big\}&\cr
&=(1-x^2)\Big\{\sqrt{a^2+1}xU_{n-2}(x)+aT_{n-1}(x)\Big\}^2&\cr
&=(1-x^2)(a^2+1)R_{n-1}^2(x).&
\end{eqnarray*}
\noindent \eop

\medskip
From the Pell identity we immediately have
\begin{Corollary} 
For $x\in[-1,1],$ 
\[|Q_n(x)|\le1\] and its maximum of 1 is attained at the endpoints $x=\pm1$ and the zeros of $R_{n-1}(x).$
\end{Corollary}

Indeed, we claim that the endpoints together with the zeros of $R_{n-1}(x)$ form the support of the optimal prediction measure. To this end we first prove that $R_{n-1}(x)$ has $n-1$ zeros in $(-1,1).$
\begin{Lemma}
The polynomials $R_n(x)$ have $n$ distinct zeros in $(-1,1)$ which interlace the extreme points of $T_n(x),$ $\cos(k\pi/n),$ $0\le k\le n.$
\end{Lemma}
\noindent {\bf Proof.} Using the fact that $T_n'(x)=nU_{n-1}(x),$ we have that at an interior
extremal point of $T_n(x),$ $\cos(k\pi/n),$ $1\le k\le (n-1),$
\begin{eqnarray*}
R_n(\cos(k\pi/n))&=&\frac{1}{\sqrt{a^2+1}} \left(\sqrt{a^2+1}zU_{n-1}(\cos(k\pi/n))+aT_{n}(\cos(k\pi/n))\right)\cr
&=&\frac{1}{\sqrt{a^2+1}} \left(0+a(-1)^k\right)\cr
&=&\frac{a}{\sqrt{a^2+1}}(-1)^k
\end{eqnarray*}
so that 
\[{\rm sgn}(R_n(\cos(k\pi/n)))=(-1)^k,\quad 1\le k\le (n-1).\]
\bigskip
Further, for $k=0,$ $\cos(k\pi/n)=1,$
\begin{eqnarray*}
R_n(1)&=&\frac{1}{\sqrt{a^2+1}} \left(\sqrt{a^2+1}U_{n-1}(1)+aT_{n}(1)\right)\cr
&=&\frac{1}{\sqrt{a^2+1}}\left(n\sqrt{a^2+1}+a\right)
\end{eqnarray*}
so that
\[{\rm sgn}(R_n(\cos(0\pi/n)))=+1=(-1)^0.\]
Similarly, for $k=n,$ $\cos(k\pi/n)=-1,$
\begin{eqnarray*}
R_n(-1)&=&\frac{1}{\sqrt{a^2+1}} \left(\sqrt{a^2+1}(-1)U_{n-1}(-1)+aT_{n}(-1)\right)\cr
&=&\frac{1}{\sqrt{a^2+1}}(n\sqrt{a^2+1}+a)(-1)^n
\end{eqnarray*}
so that also
\[{\rm sgn}(R_n(\cos(n\pi/n)))=(-1)^n.\]
The result follows. \eop

\medskip
Suppose now that $\mu_0$ is the discrete measure supported on $\pm 1$ together with the $n-1$ zeros of
$R_{n-1}(x),$ with optimal weights given by Lemma \ref{HLWeights}.

\begin{Proposition}\label{QnExtremal} The polynomials $Q_n(z)$ are of extremal growth at $z_0=ai$ relative to $K=[-1,1].$ Specifically, $Q_n(z)=-(i)^nP_n^{\mu_0,z_0}(z).$
\end{Proposition}
\noindent {\bf Proof.} Let $-1=x_0<x_1<\cdots<x_n=+1$ be the support points with corresponding Lagrange polynomials $\ell_k(z).$  We will show that
\[Q_n(x_k)=-(i)^n{\rm sgn}(\ell_k(ai)),\quad 0\le k\le n\]
using the formula
\[\ell_k(z)=\frac{\omega_n(z)}{(z-x_k)\omega_n'(x_k)},\quad \omega_n(z):=(z^2-1)R_{n-1}(z).\]
Our calculations will make use of the elementary facts that
\begin{eqnarray*}
T_n(ai)&=&\frac{(i)^n}{2}\left\{(a+\sqrt{a^2+1})^n+(a-\sqrt{a^2+1})^n\right\},\cr
U_n(ai)&=&\frac{(i)^n}{2\sqrt{a^2+1}}\left\{(a+\sqrt{a^2+1})^{n+1}-(a-\sqrt{a^2+1})^{n+1}\right\}
\end{eqnarray*}
so that
\begin{eqnarray*}
R_{n-1}(ai)&=&\frac{1}{\sqrt{a^2+1}} \left(\sqrt{a^2+1}(ai)U_{n-2}(ai)+aT_{n-1}(ai)\right)\cr
&=& (i)^{n-1}\frac{a}{\sqrt{a^2+1}}(a+\sqrt{a^2+1})^{n-1}.
\end{eqnarray*}
The endpoints are the easiest case and so we will begin with those. Specifically, for $k=0,\,x_0=-1,$
\begin{eqnarray*}
\ell_0(ai)&=&\frac{((ai)^2-1)R_{n-1}(ai)}{(ai-(-1))\omega_n'(-1)}\cr
&=&\frac{-(a^2+1)R_{n-1}(ai)}{(ai+1)(-2R_{n-1}(-1))}.
\end{eqnarray*}
Hence
\begin{eqnarray*}
{\rm sgn}(\ell_0(ai))&=&{\rm sgn}(R_{n-1}(ai))\,{\rm sgn}(R_{n-1}(-1))\,{\rm sgn}\left(\frac{1}{ai+1}\right)\cr
&=&(-i)^{n-1}(-1)^{n-1}\frac{ai+1}{\sqrt{a^2+1}}.
\end{eqnarray*}
On the other hand
\begin{eqnarray*}
Q_n(-1)&=&\frac{1}{\sqrt{a^2+1}}\left( -(a(-1)+i)T_{n-1}(-1)+\sqrt{a^2+1}(1-(-1)^2)U_{n-2}(-1)\right)\cr
&=&\frac{1}{\sqrt{a^2+1}}(a-i)(-1)^{n-1}\cr
&=&-(i)^n{\rm sgn}(\ell_0(ai)),
\end{eqnarray*}
as is easily verified.

The other endpoint $x_n=+1$ is very similar and so we suppress the details.

Consider now, $x_k,$ $1\le k\le (n-1),$ a zero of $R_{n-1}(x).$ Then
\begin{eqnarray*}
\ell_k(ai)&=&\frac{((ai)^2-1)R_{n-1}(ai)}{(ai-x_k)(x_k^2-1)R_{n-1}'(x_k)}\cr
&=&\frac{-(a^2+1)R_{n-1}(ai)}{(ai-x_k)(x_k^2-1)(R_{n-1}'(x_k))}.
\end{eqnarray*}
Hence
\begin{eqnarray*}
{\rm sgn}(\ell_k(ai))&=&{\rm sgn}(R_{n-1}(ai))\,{\rm sgn}(R_{n-1}'(x_k))\,{\rm sgn}\left(\frac{1}{ai-x_k}\right)\cr
&=&(i)^{n-1}(-1)^{k}\frac{ai-x_k}{\sqrt{a^2+x_k^2}}
\end{eqnarray*}
as ${\rm sgn}(R_{n-1}'(x_k))=(-1)^{n-1-k},$ as is easy to see.

On the other hand, from the formula for $R_{n-1}(x)$ given in Lemma \ref{Rn}, we see that
$R_{n-1}(x_k)=0$ implies that 
\[T_{n-1}(x_k)=-\frac{\sqrt{a^2+1}}{a}x_kU_{n-2}(x_k).\]
Substituting this into the formula for $Q_n$ given in Lemma \ref{Qn} we obtain
\begin{eqnarray*}
Q_n(x_k)&=&\left\{\frac{(ax_k+i)x_k}{a}+(1-x_k^2)\right\}U_{n-2}(x_k)\cr
&=& \left(\frac{a+ix_k}{a}\right) U_{n-2}(x_k).
\end{eqnarray*}
But by the Pell identity of Proposition \ref{Pell}, $|Q_n(x_k)|=1$ and so we must have
\[Q_n(x_k)=\frac{a+ix_k}{\sqrt{a^2+x_k^2}}\,{\rm sgn}\left(U_{n-2}(x_k)\right).\]
But, as the zeros of $R_{n-1}$ interlace the extreme points $T_{n-1},$ i.e., the zeros of $U_{n-2},$ it is easy to check that ${\rm sgn}(U_{n-2}(x_k)=(-1)^{n-1-k}.$ In other words,
\[Q_n(x_k)=(-1)^{n-1-k}\frac{a+ix_k}{\sqrt{a^2+x_k^2}}\]
which is easily verified to equal $-(i)^n{\rm sgn}(\ell_k(ai)),$ as claimed. \eop

\medskip
From the recursion formula for $Q_n(z)$ it is easy to see that
\[Q_n(ai)=-(i)^n\sqrt{a^2+1}(a+\sqrt{a^2+1})^{n-1}.\]
Hence we have
\begin{Proposition} For $n=1,2,\cdots$
\[\max_{p\in \C_n[z],\,\|p\|_{[-1,1]}\le 1}|p(ai)|=\sqrt{a^2+1}(|a|+\sqrt{a^2+1})^{n-1}\]
and this maximum value is attained by $Q_n(z)$ (for $a>0$).
\end{Proposition}

It is worth noting that the extremal polynomial and optimal measure, unlike the real case depend on the exterior point $z_0.$ Moreover, this extreme value is rather larger than $|T_n(ai)|.$ Indeed it is easy to show that
\[\sqrt{a^2+1}(|a|+\sqrt{a^2+1})^{n-1}-|T_n(ai)|=(\sqrt{a^2+1}-|a|)|T_{n-1}(ai)|.\]

One may of course wonder if there are similar formulas for general points $z_0\in\C\backslash[-1,1]$  (not just $z_0=ai$). However numerical experiments seem to indicate that in general there is no three-term recurrence for the extremal polynomials.

\end{document}